\documentclass{amsart}
\usepackage{graphicx, amssymb, mathrsfs}
\usepackage[all, cmtip]{xy}
\newtheorem{theorem}{Theorem}[section]

\newtheorem{corollary}[theorem]{Corollary}
\newtheorem{lemma}[theorem]{Lemma}
\newtheorem{proposition}[theorem]{Proposition}
\newtheorem{example}{Example}[section]

\theoremstyle{definition}
\newtheorem{definition}[theorem]{Definition}
\theoremstyle{remark}
\newtheorem{remark}[theorem]{Remark}

\numberwithin{equation}{section}

\newcommand{\CC}{\mathbb C}
\newcommand{\RR}{\mathbb R}

\newcommand{\fddtheta}[1]{\frac{\partial #1}{\partial \theta}}

\def\NN{{\mathbb N}}
\def\ZZ{{\mathbb Z}}

\def\PP{{\mathbb P}}

\def\cB{{\mathcal B}}
\def\cC{{\mathcal C}}

\def\cL{{\mathcal L}}

\begin{document}

\title[Generalized H\'enon Mappings and foliation by injective Brody curves]{Generalized H\'enon Mappings and foliation by injective Brody curves}%
\author{Taeyong Ahn}%
\address{(Ahn) KIAS, 85 Hoegiro, Dongdaemun-gu, Seoul 130-722, Republic of Korea}%
\email{(Ahn) tahn@kias.re.kr}
\date{\today}

\thanks{This work was supported in part by the National Research Foundation of Korea (NRF) grant funded by the Korea government (MSIP) (No. 2011-0030044).}

\begin{abstract}
We consider a finite composition of generalized H\'{e}non mappings $\mathfrak{f}:\CC^2\to\CC^2$ and its Green function $\mathfrak{g}^+:\CC^2\to\RR_{\ge 0}$ (see Section \ref{sec:Henon}). It is well known that each level set $\{\mathfrak{g}^+=c\}$ for $c>0$ is foliated by biholomorphic images of $\CC$ and each leaf is dense. In this paper, we prove that each leaf is actually an injective Brody curve in $\PP^2$ (see Section \ref{sec:Brodyness}).\
We also study the behavior of the level sets of $\mathfrak{g}^+$ near infinity.
\end{abstract}
\maketitle

\section{Introduction}

A generalized H\'enon mapping is defined by
$$
f(z, w)=(p(z)-aw, z)
$$
as a polynomial diffeomorphism of $\CC^2$ where $p:\CC\to\CC$ is a monic polynomial of degree $d\geq 2$ and $a$ is a non-zero constant. Then, we have $f^{-1}(z, w)=(w, (p(w)-z)/a)$. The class of generalized H\'enon mappings is important in the sense that according to \cite{polyauto}, finite compositions of functions of this class are the only functions that show chaotic behavior among holomorphic polynomial automorphisms of $\CC^2$.

Based on the dynamics of $f$, $\CC^2$ is divided into two sets: the set $K^+$ of points of bounded orbit and the complement $U^+$ of $K^+$. In this paper, we focus on $U^+$. For the studies of $U^+$, see Favre (\cite{Favre2000}), Hubbard and Oberste-Vorth (\cite{HO}), for instance. See also Dloussky and Oeljeklaus (\cite{DO1}, \cite{DO2}).

One of the useful methods to study the dynamics of $f$ is pluripotential theory. Let $g^+:\CC^2\to\RR_{\geq 0}$ be the Green function on $\CC^2$ for $f$ (see Section \ref{sec:Henon}). The set $U^+$ can be characterized by $U^+=\{g^+>0\}$. In \cite{HO} (see also Favre (\cite{Favre2000}), Dloussky and Oeljeklaus (\cite{DO1}, \cite{DO2})), Hubbard and Oberste-Vorth proved:
\begin{theorem}[See Theorem 7.2 in \cite{HO}]\label{thm:HO_main} The set $\cL_c:=\{g^+=c\}$ for $c>0$ is naturally foliated by Riemann surfaces. The leaves of the natural foliation of $\cL_c$ are isomorphic to $\CC$ and each leaf is dense in $\cL_c$.
\end{theorem}

The goal of this paper is to improve Theorem \ref{thm:HO_main} and find non-trivial injective Brody curves by proving the following (for injective Brody curves, see Section \ref{sec:Brodyness}):
\smallskip
\begin{theorem}\label{thm:maingeneral}
Let $\mathfrak{f}$ be a finite composition of generalized H\'enon mappings and $\mathfrak{g}^+$ its Green function. Then, the set $\{\mathfrak{g}^+=c\}$ for $c>0$ is foliated by injective Brody curves in $\PP^2$ and each leaf is dense in $\{\mathfrak{g}^+=c\}$.
\end{theorem}

For clarity, we first prove Theorem \ref{thm:BrodyLeaf_intro}, a simple version and give a sketch of the proof of Theorem \ref{thm:maingeneral} in the last section as the proofs of the two theorems are essentially the same.
\smallskip
\begin{theorem}~\label{thm:BrodyLeaf_intro}
Let $f(z, w)=(p(z)-aw, z)$ be a generalized H\'enon mapping and $g^+$ its Green function. Then, the set $\cL_c:=\{g^+=c\}$ for $c>0$ is foliated by injective Brody curves in $\PP^2$ and each leaf is dense in $\cL_c$.
\end{theorem}

References for recent developments in the space of Brody curves are Gromov(\cite{Gromov}), Eremenko(\cite{Eremenko2010}, \cite{EremenkoArxiv}), Tsukamoto (\cite{Tsukamoto2008}, \cite{Tsukamoto2009-1}, \cite{Tsukamoto2009}, \cite{Tsukamoto2012}), Matsuo and Tsukamoto (\cite{MatsuoTsukamoto}) to list a few. See also Burns and Sibony (\cite{BurnsSibony}), Do, Mai and Ninh (\cite{domaininh}).

\begin{remark}
The density property in Theorem \ref{thm:maingeneral} implies the non-triviality of the injective Brody curves.
\end{remark}
\medskip

The main idea of the proof of Theorem \ref{thm:BrodyLeaf_intro} is modifying the Brody reparametrization lemma and investigating the behavior of the set $\{g^+=c\}$ near the hyperplane at infinity. The point of our modification of the reparametrization is making the resulting limit curves of the reparametrization technique pass through a prescribed point. In some sense, we have precise information on the location of the limit curves while the Brody reparametrization lemma does not provide it in general. For the behavior of the set $\{g^+=c\}$ near the hyperplane at infinity, we prove Theorem \ref{thm:NonexistenceofHoloCurve} in Section \ref{sec:behaviour} (for $I_+$ and $K_c$, see Section \ref{sec:Henon}).
\begin{theorem}\label{thm:NonexistenceofHoloCurve}
There is no non-trivial holomorphic curve in $\PP^2$, which passes through $I_+$, and is supported in $\overline{K_c}\subseteq\PP^2$ for $c>0$.
\end{theorem}
\medskip

This paper is organized as follows: In Section 2, we briefly review generalized H\'enon mappings. In Section 3, Theorem \ref{thm:NonexistenceofHoloCurve} is proved. In Section 4, the concepts of the Brody curve and the injective Brody curve are discussed. In Section 5, we define and study a specific family of analytic discs using normal coordinates. In Section 6, Theorem \ref{thm:BrodyLeaf_intro} is proved. In Section 7, we give a sketch of the proof of Theorem \ref{thm:maingeneral}.

\subsection*{Notation}
We use $[z:w:t]$ for the homogeneous coordinate system of $\PP^2$ and $(z, w)\to[z:w:1]$ for the usual affine coordinate system of $\CC^2\subset\PP^2$ unless stated otherwise. Let $\Delta_r$ denote the disc in $\CC$ centered at the origin and of radius $r$, and $\Delta$ the unit disc in $\CC$.

We denote by $\|\cdot\|$ the standard Euclidean norm and by $ds$ the Fubini-Study metric of $\PP^2$. If necessary, we write it more precisely as $ds(p, v)$ for $p\in\PP^2$ and $v\in T_p\PP^2$. For notational simplicity, we write $\|\psi\|_{FS, \theta_0}$ to mean $ds(\psi(\theta_0), d\psi|_{\theta=\theta_0}(\frac{d}{d\theta}))$ for $\theta_0\in U$ where $U$ is an open subset of $\CC$ and $\psi:U\to\PP^2$ is a holomorphic mapping.

For a given holomorphic endomorphism $h:\CC^2\to\CC^2$, we write its $n$-th iterate as $h^n=({h^n}_1, {h^n}_2)$. By convention, $h^0(z, w)$ simply means $(z, w)$. For a given holomorphic function $P:\CC\to\CC$, $P'$ denotes the derivative of $P$.

\subsection*{Acknowledgement} The author would like to give thanks to Prof. John Erik Forn{\ae}ss and Prof. Nessim Sibony for introducing this problem and for their advice on this topic. The author thanks the referee for careful reading and for suggestions.

\section{Generalized H\'{e}non mappings}\label{sec:Henon}
Let $\PP^2$ be the $2$-dimensional complex projective space and
$$
I_+:=[0:1:0]\quad\textrm{ and }\quad I_-:=[1:0:0]
$$
in the homogeneous coordinate system of $\PP^2$. Then, $f$ has the natural extension to $\widetilde f\colon\PP^2\setminus\{I_+\}\to\PP^2\setminus\{I_+\}$ by
$$
\widetilde f([z:w:t])=\left[t^dp\left(\frac{z}{t}\right)-awt^{d-1}: zt^{d-1}: t^d\right].
$$
Similarly, $f^{-1}$ also has the natural extension to $\widetilde{f^{-1}}:\PP^2\setminus\{I_-\}\to\PP^2\setminus\{I_-\}$ by
$$
\widetilde{f^{-1}}([z:w:t])=\left[wt^{d-1}: \frac{1}{a}\left(t^dp\left(\frac{w}{t}\right)-zt^{d-1}\right): t^d\right].
$$

We recall the following notions and properties related to the dynamics of $f$ as in \cite{HO}. See also \cite{polydiffeo1}, \cite{FSHenon} and \cite{rationalPk}. Let
$$
K^\pm=\{p\in\CC^2\colon \{f^{\pm n}(p)\}\textrm{ is a bounded sequence of }n\},
$$
$K=K^+\cap K^-$, $J^\pm=\partial K^\pm$, $J=J^+\cap J^-$ and $U^\pm=\CC^2\setminus K^\pm$.
\medskip

We define the Green function on $\CC^2$ for $f$ and $f^{-1}$ by
$$
g^+(z, w):=\lim_{n\to\infty}\frac{1}{d^n}\log^+\|g^n(z, w)\|
$$
and
$$
g^-(z, w):=\lim_{n\to\infty}\frac{1}{d^n}\log^+\|g^{-n}(z, w)\|,
$$
respectively where $\log^+:=\max\{0, \log\}$. Then, $g^\pm:\CC^2\to\RR_{\geq 0}$ are non-negative, H\"older continuous on $\CC^2$, plurisubharmonic on $\CC^2$ and pluriharmonic on $U^\pm$, and satisfy
$$
g^+\circ f=dg^+\quad\textrm{ and }\quad g^-\circ f^{-1}=dg^-.
$$
It is well-known that $K^\pm=\{g^\pm=0\}$ and $U^\pm=\{g^\pm>0\}$. 
\begin{proposition}[See ~\cite{rationalPk}]\label{prop:I+I-} $K^\pm$, $U^\pm$, $I_\pm$, $\widetilde{f}$ and $\widetilde{f^{-1}}$ satisfy the following:
\begin{enumerate}
	\item $I_-$ and $I_+$ are the super-attracting fixed points of $\widetilde{f}$ and $\widetilde{f^{-1}}$, respectively,
	\item any compact subset of $U^\pm$ uniformly converges to $I_\mp$ under $\widetilde{f}$ and $\widetilde{f^{-1}}$, respectively,
	\item $\widetilde{f}(\{t=0\}\setminus I_+)=I_-$ and $\widetilde{f^{-1}}(\{t=0\}\setminus I_-)=I_+$, and
	\item $\overline{K^+}=K^+\cup I_+$ and $\overline{K^-}=K^-\cup I_-$ in $\PP^2$.
\end{enumerate}
\end{proposition}

More generally, we define
$$
K_c:=\{g^+\leq c\}, \,\,\textrm{ for }c>0 \quad\textrm{ and }\quad\cL_c:=\{g^+=c\}.
$$
\begin{proposition}[See Lemma 6.3 in ~\cite{DinhSibony}]\label{prop:Kc_closure}
$\overline{K_c}=K_c\cup \{I_+\}$.
\end{proposition}

\section{Behavior of level sets $\cL_c:=\{g^+=c\}$ near infinity}\label{sec:behaviour}
In this section, we prove Theorem \ref{thm:NonexistenceofHoloCurve}.
\medskip

\noindent
\it
Proof of Theorem \ref{thm:NonexistenceofHoloCurve}.
\rm
We prove this by contradiction. Suppose to the contrary that there exists such a non-trivial holomorphic curve $\cC$. Let $\phi:\Delta\to\PP^2$ be a parametrization of $\cC$ and we write $\phi(\theta)=[z(\theta):1:t(\theta)]$ where $z, t$ are holomorphic functions of $\theta$ with $z(0)=t(0)=0$. Notice that due to Proposition \ref{prop:Kc_closure}, $t$ cannot be identically $0$. By restricting and reparametrizing, we may further assume that $t$ vanishes only at $\theta=0$ and let $N$ denote the vanishing order of $t$ at $\theta=0$.

Since $\cC\setminus\{I_+\}=\phi(\Delta\setminus\{0\})\subset\CC^2$, $f(\cC\setminus\{I_+\})\subset\CC^2$. If $z$ is identically $0$, then the set $f(\Delta\setminus\{0\})$ accumulates at $I_-$ in $\PP^2$ near $\theta=0$. This is a contradiction to Proposition \ref{prop:Kc_closure} since $f(\cC\setminus \{I_+\})\subset f(K_c)=K_{dc}$. Hence, without loss of generality, by restricting and reparametrizing, we may assume that $z$ vanishes only at $\theta=0$ as well.

We consider limit points of $f(\Delta\setminus\{0\})$ near $\theta=0$ in $\PP^2$. We claim that $\lim_{\theta\to 0}\widetilde{f}([z(\theta): 1: t(\theta)])=I_+$. For $\theta\neq 0$, we have
$$
\widetilde{f}([z(\theta): 1: t(\theta)])=\left[p\left(\frac{z(\theta)}{t(\theta)}\right)-\frac{a}{t(\theta)}: \frac{z(\theta)}{t(\theta)}:1\right].
$$
Since $t(\theta)\to 0$ as $\theta\to 0$, points in $\CC^2$ cannot be limit points of $f(\Delta\setminus\{0\})$ near $\theta=0$. Points $f(\Delta\setminus\{0\})$ near $\theta=0$ have limit points only in the hyperplane at infinity. Since $f(\cC\setminus \{I_+\})\subset f(K_c)=K_{dc}$, Proposition \ref{prop:Kc_closure} implies our claim.

Thus, due to the Riemann removable singularity theorem, we can write
$$
f([z(\theta): 1: t(\theta)])=[z_1(\theta): 1: t_1(\theta)] \textrm{ for }\theta\neq 0,
$$
where $z_1, t_1$ are holomorphic functions of $\theta$ and vanish only at $\theta=0$. In other words, $\overline{f(\cC\setminus\{I_+\})}$ defines a non-trivial holomorphic curve in $\PP^2$, which passes through $I_+$, and is supported in $\overline{K_{dc}}\subseteq\PP^2$.

By applying this argument repeatedly, we obtain a sequence of non-trivial holomorphic curves $\{[z_n: 1: t_n]:\Delta\to \overline{K_{d^nc}}\}_{n=0, 1, 2, \cdots}$ where $z_n, t_n$ are holomorphic functions of $\theta$ and vanish only at $\theta=0$. By convention, $z_0=z$ and $t_0=t$.

Then, for $\theta$ near $0$ with $\theta\neq 0$, we have
\begin{eqnarray*}
[z_{n+1}(\theta): 1: t_{n+1}(\theta)]&=&\widetilde{f}([z_n(\theta):1:t_n(\theta)])\\
&=&\left[\frac{1}{z_n(\theta)}\right(t_n(\theta)p\left(\frac{z_n(\theta)}{t_n(\theta)}\right)-a\left): 1: \frac{t_n(\theta)}{z_n(\theta)}\right].
\end{eqnarray*}

Let $\alpha_n$ and $\beta_n$ denote the vanishing order of $z_n$ and $t_n$ at $\theta=0$. From $\lim_{\theta\to 0}t_n(\theta)/z_n(\theta)=0$, we have $\beta_n>\alpha_n$. In order to have
\begin{displaymath}
\lim_{\theta\to 0}\left[t_n(\theta)p\left(\frac{z_n(\theta)}{t_n(\theta)}\right)-a\right]=0,
\end{displaymath}
since $(z_n(\theta))^d/(t_n(\theta))^{d-1}$ is the dominating term in $t_n(\theta)p\left(z_n(\theta)/t_n(\theta)\right)$ near $\theta=0$ and $a\neq 0$, we have $d\alpha_n=(d-1)\beta_n$. Since $d$ and $d-1$ are relatively prime, $\alpha_n$ is divisible by $d-1$. From $\alpha_n<\beta_n$, we know that $\alpha_n$ is a positive integral multiple of $d-1$. Observe that $t(\theta)$ can be represented as $t(\theta)=t_{N+2}\prod_{i=0}^{N+1} z_i(\theta)$ from the relationship between $z_n$'s and $t_n$'s. It means that $t$ vanishes with order $>N$ at $\theta=0$. This is a contradiction since we assumed that $N$ is the vanishing order of $t$ at $\theta=0$. This completes the proof.
\hfill $\Box$

\section{Brody Curves}
\label{sec:Brodyness}
The Brody curve appeared in Brody's proof in \cite{Brody} that every non-Kobayashi-hyperbolic compact complex manifold contains a non-trivial holomorphic image of $\CC$. That non-trivial entire curve is called a Brody curve. It can be defined as follows:

\begin{definition}[Brody Curve]
Let $M$ be a compact complex manifold with a smooth metric $ds_M$. Let $\psi:\CC\to M$ be a non-constant holomorphic map.

The map $\psi$ is said to be \emph{Brody} if $\sup_{\theta'\in\CC}ds_M(\psi(\theta'), d\psi|_{\theta=\theta'}(\fddtheta{}))<C_\psi$ for some constant $C_\psi>0$. We call the image $\psi(\CC)$ a \emph{Brody curve} in $M$. The curve $\psi(\CC)$ is said to be \emph{injective Brody} if the parametrization $\psi$ is injective.
\end{definition}

\begin{remark}
Since the manifold is compact, a curve being Brody is independent of the choice of the metric.
\end{remark}

In the rest of the paper, for convenience, we consider the Fubini-Study metric $ds$ of $\PP^2$. With respect to the affine coordinate chart $(z, w)\in\CC^2\subset\PP^2$, the Fubini-Study metric is defined by $ds((z, w), (z', w'))=[(|z'|^2+|w'|^2+|zw'-z'w|^2)/(1+|z|^2+|w|^2)^2]^{1/2}$ for $(z, w)\in\CC^2$ and $(z', w')\in T_{(z, w)}\CC^2$.
\medskip

Below are some simple examples.
\begin{example}
Let $\alpha$ be a complex constant and $p, q$ polynomials of one complex variable $z$. Then, all curves of the form $[z: p(z): 1]$ and of the form $[p(z)\exp(z): q(z)\exp(\alpha z):1]$ are Brody in $\PP^2$. On the other hand, the mapping $z\to[\exp(z): \exp(iz^2):1]$ is not Brody.
\end{example}

\begin{example}
The map $f_n:z\to(z, \exp(z^n))$ is injective but not Brody in $\CC^2\subset\PP^2$ for $n\geq 3$. In particular, not all biholomorphic images of $\CC$ in $\PP^2$ are Brody.
\end{example}
See \cite{BurnsSibony}, \cite{domaininh} for producing a Brody curve from an entire curve. The reader is referred to the references for the space of Brody curves mentioned in the introduction.
\medskip

The following is a property of injective Brody curves. Its proof is straightforward and omitted.

\begin{proposition}
For an injective Brody curve $\cC$ in $\PP^2$, every injective parametrization of $\cC$ has its derivative bounded over $\CC$ with respect to the Fubini-Study metric. In short, a curve being injective Brody does not depend on the choice of the injective parametrization.
\end{proposition}


\section{Normal coordinates near $I_-$ and a family of analytic discs}\label{sec:analyticdiscs}
In this section, we construct a family of analytic discs using the normal coordinates introduced by Favre in \cite{Favre2000}.

\subsection{Normal coordinates}\label{subsec:normal_coordinates} We consider $\widetilde{f}$ near $I_-$. We take an affine chart $(\zeta, \omega)\to[1: \omega: \zeta]$ so that the origin maps to the point $I_-$ and let $\psi$ denote the transition map $\psi(z, w)=(1/z, w/z)$ from the usual affine coordiante chart $(z, w)\to[z: w: 1]$ to the affine coordinate chart $(\zeta, \omega)$ near $I_-$.
\medskip

In \cite{Favre2000}, Favre introduced normal coordinates (see also \cite{HO}):
\begin{proposition}[See Proposition 2.2 in \cite{Favre2000}]\label{prop:Favrecoordinate}
There exists a local biholomorphism $\phi(x, y)=(\zeta, \omega)$ near $I_-$ tangent to the identity and a holomorphic function $r$ of one complex variable vanishing at $0$ such that $\phi^{-1}\circ \widetilde{f}\circ \phi=G$ with:
$$
G(x, y)=\left(x^d, \frac{a}{d}yx^{2d-2}+x^{d-1}(1+r(x))\right).
$$
\end{proposition}
\begin{remark} Since $G$ is conjugate to $\widetilde{f}$ and the conjugation maps $\{x=0\}$ to the line at infinity, $G$ is injective outside $\{x=0\}$ in a small neighborhood of $I_-$.
\end{remark}
\begin{remark} As noticed in \cite{Favre2000}, we can extend the conjugacy $\phi$ to a global holomorphic map $\phi:(\Delta\setminus\{0\})\times\CC\to U^+$ by defining $\phi(x, y)=\psi\circ f^{-n}\circ\psi^{-1}\circ\phi\circ G^n(x, y)$ for sufficiently large $n$.
\end{remark}

In the remaining of the paper, $\phi$ denotes the global holomorphic map.

\begin{remark}
We have $\phi^{-1}\circ\psi(\{g^+=c\})=\{|x|=1/\exp(c)\}$.
\end{remark}

We choose a neighborhood of $I_-$ and constants $r_\phi, c_\phi, c_{V^+}, R$. The local biholomorphism $\phi$ is of the form $\phi(x, y)=(x(1+g(x, y)), y(1+h(x, y))$ where $g, h$ are holomorphic functions vanishing at $0$ (see p.491 of \cite{Favre2000}).
\medskip

Let $r_\phi>0$ be such that $\phi$ in Proposition \ref{prop:Favrecoordinate} is well-defined and Proposition \ref{prop:Favrecoordinate} is true for the ball $B(r_\phi, I_-)$ centered at $I_-$ and of radius $r_\phi$. Let $U_{I_-}:=\{|x|, |y|<c_\phi\}$ denote a neighborhood of $I_-$ where $c_\phi$ satisfies the following conditions:
\begin{itemize}
\item $0<c_\phi<r_\phi/4$, $c_\phi<1/100$;
\item $\sup_{U_{I_-}}\|g\|,\, \sup_{U_{I_-}}\|h\|<1/100$;
\item $c_\phi(\|g\|_{\cC^1, B(r_\phi/2, I_-)} + \|h\|_{\cC^1, B(r_\phi/2, I_-)})<1/100$.

\end{itemize}
Then, there exists a large $c_{V^+}>1/c_\phi$ and a large $R>1$ such that $\{|z|\geq R,$ $c_{V^+}|w|\leq |z|\}\subseteq \psi^{-1}\circ\phi(U_{I_-})$. By increasing $R$ if necessary as in \cite{polydiffeo1}, \cite{FSHenon} and \cite{polyauto}, we may assume that 
\begin{itemize}
\item $|z|>R$ implies either $|p(z)-aw|>c_{V^+}|z|$ or $|z|<c_{V^+}|w|$ or both;
\item $|w|>R/c_{V^+}$ implies either $|p(w)-z|/|a|>|w|$ or $c_{V^+}|w|<|z|$ or both.
\end{itemize}
By further increasing $R$, we may assume that $R$ satisfies the following approximation properties:
\begin{enumerate}
\item $|z|^d/2+c_{V^+}|z|\leq|p(z)|\leq 2|z|^d-c_{V^+}|z|$ for $|z|\geq R/c_{V^+}$;
\item $d|z|^{d-1}/2+|z|\leq|p'(z)|\leq 2d|z|^{d-1}-|z|$ for $|z|\geq R/c_{V^+}$;
\item $2^{d+5}\cdot 3^2\cdot \sqrt{41}\cdot (1+|a|)^2{c_{V^+}}^d\sqrt{1+{\|g\|_{\cC^1, U_{I_-}}}^2}\leq R$.
\end{enumerate}
For such $c_{V^+}$ and $R$, we define a filtration $\{V^+, V^-, W\}$ as in \cite{polydiffeo1}, \cite{FSHenon} and \cite{polyauto}:
\begin{eqnarray*}
&& V^+:=\{(z, w)\in\CC^2\colon R\leq|z|, c_{V^+}|w|\leq|z|\},\\
&& V^-:=\{(z, w)\in\CC^2\colon R\leq c_{V^+}|w|, |z|\leq c_{V^+}|w|\},\textrm{ and }\\
&& W:=\{(z, w)\in\CC^2\colon |z|, c_{V^+}|w|\leq R\}.
\end{eqnarray*}
Then, we have
\begin{proposition}[For example, see ~\cite{FSHenon}]\label{prop:filtration}\hspace{.1in}
$f^{\pm 1}(V^\pm)\subseteq V^\pm$ and $U^\pm=\cup_{i=0}^\infty f^{\mp i}(V^\pm)$.
\end{proposition}

\subsection{Family of analytic discs}\label{subsec:analyticdiscs}
Let $c$ be a positive constant. Notice that the level set $\{g^+=c\}$ in $V^+$ appears as $\{|x|=1/\exp(c)\}$ via $\phi^{-1}\circ\psi$ and that for any leaf $\cC$ in $\{g^+=c\}$, there exists $x_0\in\Delta\setminus \{0\}$ with $|x_0|=1/\exp(c)$ such that $\psi^{-1}\circ\phi(x_0, \theta)$ for $\theta\in\CC$ is an injective parametrization of the leaf $\cC$. Here, be careful that such $x_0$ is not uniquely determined. We prove a lemma describing that $\cC$ looks almost vertical in $V^+$ with respect to the standard coordinate chart of $\CC^2$.

\begin{lemma}\label{lem:verticality}
Let $(z(\theta), w(\theta))=\psi^{-1}\circ\phi(x_0, \theta)$.
$$
\frac{\partial z/\partial \theta}{\partial w/\partial \theta}<2\|g\|_{\cC^1, U_{I_-}} \quad\textrm{ for }|x_0|, |\theta|<c_\phi.
$$
\end{lemma}

\noindent
\it
Proof.
\rm
We can write the parametrization
$$
(z(\theta), w(\theta))=\left(\frac{1}{x_0(1+g(x_0, \theta))}, \frac{\theta(1+h(x_0, \theta))}{x_0(1+g(x_0, \theta))} \right).
$$

\begin{eqnarray*}
\frac{\partial z/\partial \theta}{\partial w/\partial \theta}&=&\frac{g'/(x_0(1+g)^2)}{(1+h)/(x_0(1+g))-\theta(h'(1+g)-g'(1+h))/(x_0(1+g)^2)}\\
&=&\frac{g'}{(1+h)(1+g)-\theta(h'(1+g)-g'(1+h))}
\end{eqnarray*}
where $g'=dg(x_0,\theta)/d\theta$ and $h'=dh(x_0,\theta)/d\theta$. From our choice of $c_\phi$, the last quantity is bounded by $2\|g\|_{\cC^1, U_{I_-}}$ for $|\theta|<c_\phi$.
\hfill $\Box$

As we will use the constant frequently, we denote by $c_g:=2\|g\|_{\cC^1, U_{I_-}}$.
\medskip

Now, we define a family of analytic discs. In the remaining of this section, we assume that $c$ is sufficiently large so that
\begin{itemize}
\item $c>\max_W g^+$;
\item $(|a|c_\phi/d+\max_{U_{I_-}}|r| +1)<c_\phi\exp(c)/8$ and $R<\exp(c)$
\end{itemize}
where $r$ is as in Proposition \ref{prop:Favrecoordinate}.

We consider a leaf $\cC\subset\{g^+=c\}$. Then, as just discussed, we can find $s\in\CC$ with $|s|=1/\exp(c)$ such that $\psi^{-1}\circ\phi(s, \CC)=\cC$. For each $n\in\NN$, we define an analytic disc $\phi_{s,n}:\Delta\to\cC\subset \CC^2$ by
$$
\phi_{s, n}(\theta)=f^{-n}\circ\psi^{-1}\circ\phi(s^{d^n}, c_\phi\theta/2+\theta_{s,n})
$$
where $c_\phi$ is the constant in the definition of the neighborhood $U_{I_-}$ and $\theta_{s,n}:={G^n}_2(s, 0)$. These analytic discs are well defined from our choice of $c$. Then, for all $n$, $\phi_{s,n}(0)=\phi_{s,0}(0)\in\cC$ and we denote this point by $P$. We write $\Phi_{s,n}:=\phi_{s,n}(\Delta)$. It is not difficult to see that $\Phi_{s,n}\subset\Phi_{s, n+1}$. The sequence $\{\phi_{s,n}\}$ of analytic discs is the desired family of analytic discs and will be used for the proof of Theorem \ref{thm:BrodyLeaf_intro}.
\medskip

We investigate the Fubini-Study metric of the map $\phi_{s,n}$ over $\Delta$. As $\phi$ is a local biholomorphism tangent to the identity, it suffices to consider the action of $f^{-1}$ on the Fubini-Study metric of $f^{-i}\circ\psi^{-1}\circ\phi(s^{d^n}, c_\phi\theta/2+\theta_{s,n})$ (or equivalently $f^{n-i}\circ \phi_{s,n}$) for $i=0, \ldots, n-1$. For notational convenience, we define the following sets: for $0\leq i\leq n$,
$$
\Theta^s_{n,i}:=\{\theta\in\Delta\colon f^j(\phi_{s,n}(\theta))\in V^+ \textrm{ for all }j\textrm{ with }i\leq j\leq n\}.
$$
Then, since $f(V^+)\subseteq V^+$ from Proposition \ref{prop:filtration}, $\Theta^s_{n,i}\subseteq\Theta^s_{n,i+1}$ is clear. In particular, $\Theta^s_{n,0}$ is the set of $\theta$'s with $\phi_{s,n}(\theta)\in V^+$.

Due to Proposition \ref{prop:filtration}, $g^+\circ f=dg^+$ and our choice of $c$, we know that $f^i(\phi_{s,n}(\Delta))$ lies in $V^+\cup V^-$ and it is enough to consider $f^{-1}$ in three cases : 
\begin{enumerate}
\item[{\bf Case i}]$f^{-1}:f^i(\phi_{s,n}(\Theta^s_{n,i-1}))\to f^{i-1}(\phi_{s,n}(\Delta))\cap V^+$ maps points in $f^i(\phi_{s,n}(\Delta))$ from $V^+$ to $V^+$,
\item[{\bf Case ii}]$\displaystyle f^{-1}:f^i(\phi_{s,n}(\Theta^s_{n,i}\setminus\Theta^s_{n,i-1}))\to f^{i-1}(\phi_{s,n}(\Delta))\cap V^-$ maps points in $f^i(\phi_{s,n}(\Delta))$ from $V^+$ to $V^-$, and
\item[{\bf Case iii}]$\displaystyle f^{-1}:f^i(\phi_{s,n}(\Delta\setminus\Theta^s_{n,i}))\to f^{i-1}(\phi_{s,n}(\Delta))\cap V^-$ maps points in $f^i(\phi_{s,n}(\Delta))$ from $V^-$ to $V^-$.
\end{enumerate}

We study Case i. We start by proving the following lemma.
\begin{lemma}\label{lem:V+} Let $i, n\in\ZZ$ be given such that $0\leq i\leq n$. Let $(z, w)=f^i\circ\phi_{s,n}$ and $(z', w')=[d(f^i\circ\phi_{s,n})](\frac{d}{d\theta})$. Then, for any $\theta\in\Theta^s_{n,i}$, we have $|z'|\leq c_g|w'|$ and
\begin{eqnarray}
&&\frac{|w'|^2}{4(1+|z|^2+|w|^2)}\leq\frac{|z'|^2+|w'|^2+|z'w-zw'|^2}{(1+|z|^2+|w|^2)^2}\leq\frac{2(1+{c_g}^2)|w'|^2}{1+|z|^2+|w|^2}.
\end{eqnarray}
\end{lemma}

\noindent
\it
Proof.
\rm
Lemma \ref{lem:verticality} implies the first assertion. Since $c_{V^+}|w|\leq|z|$ and $c_{V^+}>1/c_\phi$,
\begin{eqnarray}
\nonumber && \frac{|z'|^2+|w'|^2+|z'w-zw'|^2}{(1+|z|^2+|w|^2)^2}\geq\frac{|w'|^2+|zw'|^2/2}{(1+|z|^2+|w|^2)^2}\geq\frac{|w'|^2}{4(1+|z|^2+|w|^2)}.
\end{eqnarray}
On the other hand, 
\begin{eqnarray}
\nonumber && \frac{|z'|^2+|w'|^2+|z'w-zw'|^2}{(1+|z|^2+|w|^2)^2}\leq\frac{|z'|^2+|w'|^2+2(|z'w|^2+|zw'|^2)}{(1+|z|^2+|w|^2)^2}\\
\nonumber && \leq2\frac{|z'|^2+|w'|^2}{1+|z|^2+|w|^2}\leq\frac{2(1+{c_g}^2)|w'|^2}{1+|z|^2+|w|^2}.
\end{eqnarray}
\hfill $\Box$
\medskip

The following proposition implies that in Case i, $f^{-1}$ increases the Fubini-Study metric at least by a fixed ratio. 
\begin{lemma}\label{prop:orderV+}
Let $1\leq i\leq n$.
$$
1<\frac{d^2}{96(1+{c_g}^2)\cdot 2^{2d-2}|a|^2}|s|^{-(4-4/d)d^i}\leq\inf_{\theta\in\Theta^s_{n,i-1}}\frac{\|f^{i-1}\circ\phi_{s,n}\|^2_{FS, \theta}}{\|f^i\circ\phi_{s,n}\|^2_{FS, \theta}}.
$$
\end{lemma}

\noindent
\it
Proof.
\rm
Let $\theta\in\Theta^s_{n,i-1}$. For simplicity, write $(z, w)=f^i\circ\phi_{s,n}$, $(z', w')=[d(f^i\circ\phi_{s,n})](\frac{d}{d\theta})$, and $(z_*, w_*)=f^{i-1}\circ\phi_{s,n}$ at $\theta$ with respect to the usual coordinate system of $\CC^2$.

By Lemma \ref{lem:V+}, we have
\begin{eqnarray*}
\|f^i\circ\phi_{s,n}\|^2_{FS, \theta} \leq \frac{2(1+{c_g}^2)|w'|^2}{1+|z|^2+|w|^2}\leq\frac{2(1+{c_g}^2)|w'|^2}{|z|^2}.
\end{eqnarray*}

On the other hand, by Lemma ~\ref{lem:V+} with $(z,w),(w, (p(w)-z)/a)\in V^+$ and the chain rule,
\begin{flalign}
\nonumber\|f^{i-1}\circ\phi_{s,n}\|^2_{FS, \theta}&\geq \frac{|-z'/a+p'(w)w'/a|^2}{4(1+|w|^2+|(p(w)-z)/a|^2)}\\
\nonumber&\geq \frac{|-z'/a+p'(w)w'/a|^2}{12|w|^2}\geq\frac{(|p'(w)|-c_g)^2|w'|^2}{12|a|^2|w|^2}.
\end{flalign}
Hence,
$$
\inf_{\theta\in\Theta^s_{n,i-1}}\frac{\|f^{i-1}\circ\phi_{s,n}\|^2_{FS, \theta}}{\|f^i\circ\phi_{s,n}\|^2_{FS, \theta}}\geq\frac{(|p'(w)|-c_g)^2|z|^2}{24(1+{c_g}^2)|a|^2|w|^2}.
$$

Notice that $\phi$ is a local biholomorphism on $U_{I_-}$ and $V^+\subseteq\psi^{-1}\circ\phi(U_{I_-})$. For $\theta\in\Theta^s_{n,i-1}$, $G^{i-1-n}(s^{d^n}, \theta)$ is uniquely determined. So, we can use the conjugacy relationship. Proposition 2.3 in \cite{Favre2000} implies
$$
f^{i-1}\circ\phi_{s,n}=f^{i-1-n}\circ\psi^{-1}\circ\phi(s^{d^n}, \theta)=\psi^{-1}\circ\phi\circ G^{i-1-n}(s^{d^n}, \theta)=\psi^{-1}\circ\phi(s^{d^{i-1}}, q(\theta))
$$
where $q$ is a function of $\theta$ such that $(s^{d^{i-1}}, q(\theta))\in U_{I_-}$ for $\theta\in \Theta^s_{n,i-1}$.

By our choice of $U_{I_-}$, we have that
$$
\left(1-\frac{1}{100}\right)s^{d^{i-1}}\leq \frac{1}{|z_*|}=\frac{1}{|w|}\leq \left(1+\frac{1}{100}\right)s^{d^{i-1}}.
$$

In the same way, we have
\begin{equation}\label{ineq:range_z}
\left(1-\frac{1}{100}\right)s^{d^{i}}\leq \frac{1}{|z|}\leq \left(1+\frac{1}{100}\right)s^{d^{i}}.
\end{equation}

Therefore, from our choice of $R$, we have
\begin{eqnarray*}
\frac{(|p'(w)|-c_g)^2|z|^2}{24(1+{c_g}^2)|a|^2|w|^2}&\geq&\frac{|d|w|^{d-1}/2|^2|z|^2}{24(1+{c_g}^2)|a|^2|w|^2}=\frac{d^2}{96(1+{c_g}^2)|a|^2}|w|^{2d-4}|z|^2\\
&\geq&\frac{d^2}{96(1+{c_g}^2)\cdot 2^{2d-2}|a|^2}|s|^{-(4-4/d)d^i}>1.
\end{eqnarray*}
\hfill $\Box$
\medskip

We consider Case iii.
\begin{lemma}\label{prop:V-toV-} Let $1\leq i\leq n$.
$$
\sup_{\theta\in\Delta\setminus\Theta^s_{n,i}}\frac{\|f^{i-1}\circ\phi_{s,n}\|^2_{FS, \theta}}{\|f^i\circ\phi_{s,n}\|^2_{FS, \theta}}\leq \frac{2^4\cdot 3^2\cdot 41 d^2|a|^2{c_{V^+}}^{2d}}{R^{2d-4}}.
$$
Here, this upper bound is bounded above by the lower bound in Lemma \ref{prop:orderV+}.
\end{lemma}

\noindent
\it
Proof.
\rm
Let $\theta\in\Delta\setminus\Theta^s_{n,i}$. For simplicity, write $(z, w)=f^i\circ\phi_{s,n}$, $(z', w')=[d(f^i\circ\phi_{s,n})](\frac{d}{d\theta})$ at $\theta$ with respect to the usual coordinate system of $\CC^2$. Then, $|z|\leq c_{V^+}|w|$ and $R\leq c_{V^+}|w|$.

We have
$$
\|f^i\circ\phi_{s,n}\|^2_{FS, \theta}=\frac{|z'|^2+|w'|^2+|z'w-zw'|^2}{(1+|z|^2+|w|^2)^2}\geq\frac{|z'|^2+|w'|^2}{({3c_{V^+}}^2|w|^2)^2}.
$$

On the other hand, from our choice of $R$, we have
\begin{align*}
&\|f^{i-1}\circ\phi_{s,n}\|^2_{FS, \theta}\\
&=\frac{|w'|^2+|-z'/a+p'(w)w'/a|^2+|w(-z'/a+p'(w)w'/a)-(p(w)-z)w'/a|^2}{(1+|w|^2+|(p(w)-z)/a|^2)^2}\\
&\leq|a|^2\frac{|aw'|^2+(|z'|+|p'(w)w'|)^2+[|wz'|+(|p'(w)w|+|p(w)|+|z|)|w'|]^2}{|p(w)-z|^4}\\
&\leq|a|^2\frac{|a|^2+(1+|p'(w)|)^2+(|w|+|p'(w)w|+|p(w)|+|z|)^2}{|p(w)-z|^4}\max\{|z'|, |w'|\}^2\\
&\leq|a|^2\frac{|a|^2+(2d|w|^{d-1})^2+(|w|+2d|w|^d+2|w|^d+c_{V^+}|w|)^2}{|w|^{4d}/2^4}\max\{|z'|, |w'|\}^2\\
& \leq \frac{2^4\cdot 41 d^2|a|^2\max\{|z'|, |w'|\}^2}{|w|^{2d}}.
\end{align*}
Thus, the ratio is bounded by $2^4\cdot 3^2\cdot 41 d^2|a|^2{c_{V^+}}^{2d}/R^{2d-4}$.
From our choice of $R$, this is bounded above by the lower bound of Lemma \ref{prop:orderV+}.
\hfill $\Box$
\medskip

We finally consider Case ii. Here, we show that the increasing rate of the Fubini-Study metric in this case is dominated by the increasing rate of the Fubini-Study metric in Case i up to a uniform constant.

\begin{lemma}\label{prop:V+toV-}
For all $i, n$ such that $1\leq i\leq n$, there exists a constant $C_s>0$ independent of $i, n$ such that
$$
\sup_{\theta\in\Theta^s_{n,i}\setminus\Theta^s_{n,i-1}}\frac{\|f^{i-1}\circ\phi_{s,n}\|^2_{FS, \theta}}{\|f^i\circ\phi_{s,n}\|^2_{FS, \theta}}\leq C_s\inf_{\theta\in\Theta^s_{n,i-1}}\frac{\|f^{i-1}\circ\phi_{s,n}\|^2_{FS, \theta}}{\|f^i\circ\phi_{s,n}\|^2_{FS, \theta}}.
$$
\end{lemma}

\noindent
\it
Proof.
\rm
Let $\theta\in\Theta^s_{n,i}\setminus\Theta^s_{n,i-1}$. Then, $f^i(\phi_{s,n}(\theta))\in V^+$ but $f^{i-1}(\phi_{s,n}(\theta))\in V^-$. For simplicity, we write $(z, w)=f^i\circ\phi_{s,n}$, $(z', w')=[d(f^i\circ\phi_{s,n})](\frac{d}{d\theta})$, $(z_*, w_*)=f^{i-1}\circ\phi_{s,n}$ and $(z_*', w_*')=[d(f^{i-1}\circ\phi_{s,n})](\frac{d}{d\theta})$ at $\theta$ with respect to the usual coordinate system of $\CC^2$. Then, $|z|\geq R$, $c_{V^+}|w|\leq |z|$, $|z'|\leq c_g|w'|$ and $|w|\leq c_{V^+}|(p(w)-z)/a|$.

From Lemma ~\ref{lem:V+}, we have
$$
\|f^i\circ\phi_{s,n}\|^2_{FS, \theta}\geq\frac{|w'|^2}{4(1+|z|^2+|w|^2)}\geq\frac{|w'|^2}{12|z|^2}.
$$
We have
\noindent$\displaystyle\|f^{i-1}\circ\phi_{s,n}\|^2_{FS,\theta}$
\begin{align*}
&=\frac{|w'|^2+|-z'/a+p'(w)w'/a|^2+|(p(w)-z)w'/a-w(-z'/a+p'(w)w'/a)|^2}{(1+|w|^2+|(p(w)-z)/a|^2)^2}\\
&\leq\frac{|w'|^2+(|z'/a|+|p'(w)w'/a|)^2+(|(p(w)-z)w'/a-wp'(w)w'/a|+|wz'/a|)^2}{(1+|w|^2+|(p(w)-z)/a|^2)^2}\\
&\leq\frac{|w'|^2+(c_g|w'/a|+|p'(w)w'/a|)^2+(|(p(w)-z)w'/a-wp'(w)w'/a|+c_g|ww'/a|)^2}{(1+|w|^2+|(p(w)-z)/a|^2)^2}\\
&\leq\frac{(1+(c_g/|a|+|p'(w)/a|)^2+2|(p(w)-z)/a|^2+2(c_g+|p'(w)|)^2|w/a|^2)|w'|^2}{(1+|w|^2+|(p(w)-z)/a|^2)^2}\\
&\leq\frac{(2+2|(p(w)-z)/a|^2+(c_g+|p'(w)|)^2(1/|a|^2+2|w/a|^2))|w'|^2}{(1+|w|^2+|(p(w)-z)/a|^2)^2}\\
&\leq\frac{2|w'|^2}{1+|w|^2+|(p(w)-z)/a|^2}\\
&\hspace{1in}+\frac{(-2|w|^2+(c_g+|p'(w)|)^2(1/|a|^2+2|w/a|^2))|w'|^2}{(1+|w|^2+|(p(w)-z)/a|^2)^2}.
\end{align*}

Thus, we are looking for an upper bound of the following:
\begin{eqnarray}\label{qty:RatioV+V-}
&&\frac{24|z|^2}{1+|w|^2+|(p(w)-z)/a|^2}\\
&&\nonumber\hspace{1in}+\frac{12(-2|w|^2+(c_g+|p'(w)|)^2(1/|a|^2+2|w/a|^2))|z|^2}{(1+|w|^2+|(p(w)-z)/a|^2)^2}.
\end{eqnarray}
\medskip

We estimate upper bounds of \eqref{qty:RatioV+V-} in 3 cases:
\begin{enumerate}
	\item[\textbf{Case 1}] $|w|<R$,
	\item[\textbf{Case 2}] $R\leq|w|< 8^{1/d}|s|^{-d^{i-1}}$ and
	\item[\textbf{Case 3}] $8^{1/d}|s|^{-d^{i-1}}\leq|w|$.
\end{enumerate}
\medskip
%
\noindent\textbf{Case 1.} As $(z, w)\in V^+$, from \eqref{ineq:range_z}, it is easy to see that
$$
\eqref{qty:RatioV+V-}\leq (48+M_1)|s|^{-2d^i}
$$
where $M_1=\max_{|w|\leq R}12(-2|w|^2+(c_g+|p'(w)|)^2(1/|a|^2+2|w/a|^2)).$
\medskip

%
%

\noindent\textbf{Case 2.} In this case, from our choice of $R$, we have
$$
\frac{1}{2}d|w|^{d-1}\leq|p'(w)|\leq 2d|w|^{d-1}.
$$
Again, as $(z, w)\in V^+$, from \eqref{ineq:range_z} and from our choice of $R$, we have
\begin{eqnarray*}
\eqref{qty:RatioV+V-}&\leq&24|z|^2+M_{a, c_g, R}\frac{|p'(w)|^2|w|^2|z|^2}{|w|^4}\leq 24|z|^2+4d^2M_{a, c_g, R}|w|^{2d-4}|z|^2\\
&\leq&48|s|^{-2d^i}+\widetilde{M}_{a, c_g, R}|s|^{-(4-4/d)d^i}
\end{eqnarray*}
where $M_{a, c_g, R}$ and $\widetilde{M}_{a, c_g, R}$ are constants depending only on $a, c_g$ and $R$.
\medskip

\noindent\textbf{Case 3.} In this case, with \eqref{ineq:range_z}, we have
$$
\frac{1}{2}d|w|^{d-1}\leq|p'(w)|\leq 2d|w|^{d-1}\quad\textrm{ and }\quad\frac{1}{4}|w|^d-|z|\geq 0
$$
and therefore,
$$
|p(w)-z|\geq|p(w)|-|z|\geq \frac{1}{2}|w|^d-|z|\geq \frac{1}{4}|w|^d+\left(\frac{1}{4}|w|^d-|z|\right)\geq \frac{1}{4}|w|^d.
$$
As in Case 2, from \eqref{ineq:range_z}, we have
\begin{eqnarray*}
\eqref{qty:RatioV+V-}&\leq&24|z|^2+M_{a, c_g, R}\frac{|p'(w)|^2|w|^2|z|^2}{|p(w)-z|^4}\leq 24|z|^2+ 4^4\cdot 4d^2M_{a, c_g, R}\frac{|z|^2}{|w|^{2d}}\\
&\leq& 48|s|^{-2d^i}+\widetilde{\widetilde{M}}_{a, c_g, R}
\end{eqnarray*}
where $M_{a, c_g, R}$ and $\widetilde{\widetilde{M}}_{a, c_g, R}$ are constants depending only on $a, c_g$ and $R$.
\medskip

Since the upper bounds of \eqref{qty:RatioV+V-} in the above 3 cases have the same or less order of $|s|^{-1}$ than the lower bound in Lemma \ref{prop:orderV+}, we can find the maximum ratio of the bounds of \eqref{qty:RatioV+V-} to the lower bound in Lemma \ref{prop:orderV+}. This maximum ratio is the desired constant for $C_s$. Independence is clear from the proof.
\hfill $\Box$
\bigskip

Since $f^{-1}(V^-)\subseteq V^-$, Lemma \ref{prop:orderV+}, Lemma \ref{prop:V-toV-} and Lemma \ref{prop:V+toV-} together with $\Theta^s_{n,n}=\Delta$ imply that the supremum of the Fubini-Study metric on $\Delta\setminus \Theta^s_{n,0}$ is dominated by the infimum of the Fubini-Study metric on $\Theta^s_{n,0}$ when $n\in\NN$. Hence, we are going to compare the values of the Fubini-Study metric over $\Theta^s_{n,0}$.

\begin{lemma}\label{prop:V+toV+}
Let $0\leq i\leq n$. Then, there exists a constant $C_{s, V^+}>0$ independent of $i$ and $n$ such that
\begin{displaymath}
\frac{\sup_{\theta\in\Theta^s_{n,i}}\|f^i\circ\phi_{s,n}\|^2_{FS, \theta}}{\inf_{\theta\in\Theta^s_{n,i}}\|f^i\circ\phi_{s,n}\|^2_{FS, \theta}}<C_{s, V^+}.
\end{displaymath}
Indeed, the constant is determined by the local biholomorphism $\phi$ near the origin.
\end{lemma}

\noindent
\it
Proof.
\rm
We denote $(z(\theta), w(\theta))=f^i\circ\phi_{s,n}(\theta)$. Since $f^i\circ\phi_{s,n}(\theta)=f^{i-n}\circ\psi^{-1}\circ\phi(s^{d^n}, c_\phi\theta/2+\theta_{s,n})\in V^+$ for $\theta\in\Theta^s_{n,i}$, $G^{i-n}(s^{d^n}, c_\phi\theta/2+\theta_{s,n})$ is well defined and so we have
$$
f^i\circ\phi_{s,n}(\theta)=f^{i-n}\circ \psi^{-1}\circ\phi(s^{d^n}, c_\phi\theta/2+\theta_{s,n})=\psi^{-1}\circ\phi\circ G^{i-n}(s^{d^n}, c_\phi\theta/2+\theta_{s,n}).
$$
From Proposition 2.3 in \cite{Favre2000}, we have
$$
G^{i-n}(s^{d^n}, c_\phi\theta/2+\theta_{s,n})=\left(s^{d^i}, \frac{c_\phi\theta/2+\theta_{s,n}-Q_{n-i}(s^{d^i})}{(a/d)^{n-i}s^{q_{n-i}}}\right)
$$
where $q_{n-i}=2(d^{n-i}-1)$ and $Q_{n-i}$ is a polynomial of one complex variable. We obtain $(z'(\theta), w'(\theta))=d(\psi^{-1}\circ\phi)|_{G^{i-n}(s^{d^n}, c_\phi\theta/2+\theta_{s,n})}(0, c_\phi d^{n-i}/(2a^{n-i} s^{q_{n-i}}))$.

From Proposition \ref{prop:Favrecoordinate}, we have
$$
\psi^{-1}\circ\phi(x,y)=\left(\frac{1}{x(1+g(x, y))}, \frac{y(1+h(x, y))}{x(1+g(x, y))}\right),
$$ and therefore, from our choice of $c_\phi$, we have for $\theta\in\Theta^s_{n,i}$,
$$
\left(\frac{99}{101|s|^{d^i}}-\frac{1}{2|s|^{d^i}}\right)\frac{c_\phi d^{n-i}}{2|a|^{n-i}|s|^{q_{n-i}}}< |w'(\theta)|<\left(\frac{101}{99|s|^{d^i}}+\frac{1}{2|s|^{d^i}}\right)\frac{c_\phi d^{n-i}}{2|a|^{n-i}|s|^{q_{n-i}}}
$$
Note that the left-side of the inequality is not $0$.
As in \eqref{ineq:range_z} of Lemma \ref{prop:orderV+}, we have
$$
\left(1-\frac{1}{100}\right)s^{d^{i}}\leq \frac{1}{|z(\theta)|}\leq \left(1+\frac{1}{100}\right)s^{d^{i}}.
$$
Together with the above two inequalities, Lemma \ref{lem:V+} and $(z(\theta), w(\theta))\in V^+$ complete the proof. 
\hfill$\Box$
\medskip

Summarizing the results from Lemma \ref{prop:orderV+}, Lemma \ref{prop:V-toV-}, Lemma \ref{prop:V+toV-} and Lemma \ref{prop:V+toV+}, we obtain the following lemma for the proof of Theorem \ref{thm:BrodyLeaf_intro}:
\begin{lemma}\label{lem:mainBrody}
The sequence $\displaystyle \left\{\frac{\sup_{\theta\in\Delta}\|\phi_{s,n}\|_{FS, \theta}}{\|\phi_{s,n}\|_{FS, 0}}\right\}$ as a sequence of $n$ is bounded above by a constant $M_s>0$.
\end{lemma}

We close this section by proving the following lemma for the next section.
\begin{lemma}\label{lem:divergenceFS}
We have $\|\phi_{s,n}\|_{FS, 0}\to\infty$ as $n\to\infty$. 
\end{lemma}

\noindent
\it
Proof.
\rm
Note that for all $n$, $\phi_{s,n}(0)=P\in V^+$ and $G^{-n}(s^{d^n}, \theta_{s,n})=(s, 0)$ in $U_{I_-}$. So, as in Lemma \ref{prop:V+toV+}, we have
$$
\phi_{s,n}(\theta)=f^{-n}\circ\psi^{-1}\circ \phi(s^{d^n}, c_\phi\theta/2+\theta_{s,n})=\psi^{-1}\circ\phi\circ G^{-n}(s^{d^n}, c_\phi\theta/2+\theta_{s,n}).
$$
and 
$$
G^{-n}(s^{d^n}, c_\phi\theta/2+\theta_{s,n})=\left(s, \frac{c_\phi\theta/2+\theta_{s,n}-Q_n(s)}{(a/d)^ns^{q_n}}\right)
$$
where $q_n=2(d^n-1)$ and $Q_n$ is a polynomial of one complex variable.

We obtain $\partial \phi_{s,n}/\partial\theta=d(\psi^{-1}\circ\phi)|_{(s, 0)}(0, c_\phi d^n/(2a^n s^{q_n}))$ at $\theta=0$. Hence, from our choice of $c$ in $|s|=1/\exp(c)$, $|\partial \phi_{s,n}/\partial\theta|\to\infty$ as $n\to\infty$ at $\theta=0$. Lemma \ref{lem:V+} completes the proof.
\hfill $\Box$

\section{Proof of Theorem \ref{thm:BrodyLeaf_intro}}\label{sec:proof}
In this section, we prove Theorem \ref{thm:BrodyLeaf_intro}. We first prove the following lemma:
\begin{lemma}\label{lem:sitinleaf}
Let $\xi:\CC\to\cL_c$ be a holomorphic mapping. Then, the image $\xi(\CC)$ should sit inside a single leaf of the foliation of $\cL_c$.
\end{lemma}

\noindent
\it
Proof.
\rm
Suppose the contrary. Then, there exists an open subset $U\subset\CC$ with compact closure such that $\xi(U)$ is not contained in a single leaf. Since $\xi(U)$ has compact closure and $\xi(U)\subset\cL_c\subset U^+$, we can find a sufficiently large number $n$ such that $f^n(\xi(U))\subset V^+$. Since $f^n(\xi(U))$ does not still lie inside a single leaf, ${\phi^{-1}}_1\circ f^n\circ \xi:U\to\CC$ is a non-constant holomorphic function over the open subset $U$. So, the open mapping theorem of one complex variable implies that ${\phi^{-1}}_1\circ f^n\circ \xi(U)$ should be open in $\CC$. However, $|{\phi^{-1}}_1\circ f^n\circ \xi|\equiv 1/\exp(d^nc)$ on $U$. This is a contradiction. So, the lemma is proved.
\hfill $\Box$
\bigskip

\noindent
\it
Proof of Theorem \ref{thm:BrodyLeaf_intro}.
\rm
First, assume that $c$ is sufficiently large as in Subsection \ref{subsec:analyticdiscs} and let $\cC$ be a leaf in the natural foliation of the level set $\cL_c:=\{g^+=c\}$. The general case will be discussed later.
\medskip

We first modify the Brody reparametrization lemma to prove that there exists a Brody curve $\cB_s$ in $\cC$. Recall the family $\{\phi_{s,n}:\Delta\to\cC\}_{n=0}^\infty$ of analytic discs in Subsection \ref{subsec:analyticdiscs} and $\phi_{s,n}(0)=P\in \cC$ for all $n$. Let $R_n=\|\phi_{s,n}\|_{FS,0}$ and define a new sequence $\{k_n(\theta):=\phi_{s,n}(\theta/R_n)\}$. Then Lemma \ref{lem:mainBrody} implies that over $\Delta_{R_n/2}$,
$$
\|k_n\|_{FS, \theta}=\frac{\|\phi_{s,n}\|_{FS, \theta/R_n}}{R_n}\leq M_s.
$$
The set $\{k_n'\}$ of derivatives is uniformly bounded with respect to the Fubini-Study metric and $k_n(0)=P$ for every $n$. Lemma \ref{lem:divergenceFS} implies that $R_n\to\infty$. So, we can apply the normal family argument to $\{k_n(\theta)\}$ to obtain limit maps from $\CC$ to $\PP^2$ as in the Brody reparametrization lemma. Observe that $\|k_n\|_{FS, 0}=1$. Thus, the limit maps of $\{k_n(\theta)\}$ are Brody maps and pass through $P\in\cC$ at the origin.

We denote a limit map by $\kappa_s:\CC\to\overline{\cC}=\cL_c\cup\{I_+\}$ and its image by $\cB_s$. 
By Theorem \ref{thm:NonexistenceofHoloCurve}, we have $\cB_s\subseteq\cL_c$. Lemma \ref{lem:sitinleaf} implies that $\cB_s$ should sit in one single leaf of the natural foliation of $\cL_c$. Then, $\kappa_s(0)=P\in\cC$ means that $\cB_s$ should sit inside $\cC$.
\medskip

Next, we prove the injectivity of the limit map $\kappa_s:\CC\to\cC$. Without loss of generality, by passing to a convergent subsequence, we may assume that $\{k_n\}$ is a locally uniformly convergent sequence to $\kappa_s$. We prove that for $\theta_a, \theta_b\in\CC$ with $\theta_a\neq\theta_b$, $\kappa_s(\theta_a)\neq\kappa_s(\theta_b)$.

Let $r>0$ be such that $\theta_a, \theta_b\in\Delta_r\subset\CC$. Let $U_{\Delta_r}$ denote an open subset of $U^+$ with compact closure such that $\overline{\kappa_s(\Delta_r)}\subset U_{\Delta_r}$. By Proposition \ref{prop:I+I-}, there exists a large number $N\in\mathbb{N}$ such that $f^N(\overline{U_{\Delta_r}})\subset V^+$. 
Since $k_n\to\kappa_s$ locally uniformly, we can find another large number $N'\in\mathbb{N}$ such that for all $n\geq N'$, $k_n\to\kappa_s$ is uniform over $\Delta_r$ and $\overline{k_n(\Delta_r)}\subset U_{\Delta_r}$. From our construction, it is clear that for all $n\geq N'$, $k_n(\Delta_{R_n})\subset \cC$. Then, ${\phi^{-1}}_1\circ\psi\circ f^N\circ k_n=s^{d^N}$ and ${\phi^{-1}}_1\circ\psi\circ f^N\circ \kappa_s=s^{d^N}$ on $\Delta_r$. Then, as $k_n$ is an injective parametrization, ${\phi^{-1}}_2\circ\psi\circ f^N\circ k_n$ is injective on $\Delta_r$ for all $n\geq N'$. Since $k_n\to\kappa_s$ uniformly on $\Delta_r$, by the Hurwitz theorem of one complex variable, ${\phi^{-1}}_2\circ\psi\circ f^N\circ \kappa_s$ is either injective or constant. Since the limit map $\kappa_s$ is not constant, we have ${\phi^{-1}}_2\circ\psi\circ f^N\circ \kappa_s$ is injective and so, $\kappa_s(\theta_a)\neq \kappa_s(\theta_b)$.
\medskip

Finally, we show that $\cB_s=\cC$. So far, we have proved that $\cB_s$ is a biholomorphic image of $\CC$ in $\cC$. Suppose to the contrary that $\cB_s\neq\cC$. Recall from Theorem \ref{thm:HO_main} that $\cC$ is a biholomorphic image of $\CC$. Then, we can find a biholomorphic mapping from $\CC$ onto a proper subset of $\CC$. This is a contradiction to the uniformization theorem of one complex variable. Hence, $\cC$ itself is an injective Brody curve $\cB_s$.
\medskip

The general case is obtained by making $c$ large enough by applying $f$ sufficiently many times so that we can apply the above arguments. Indeed, it is not difficult to see that if a holomorphic map $\psi:\CC \to \CC^2\subset\PP^2$ is injective Brody, then $f^{-1}\circ \psi$ is also injective Brody.
\hfill $\Box$

\begin{remark}
In ~\cite{Fornaess}, short $\CC^k$'s were first introduced. A domain is called a \emph{short $\CC^k$} if it an increasing union of holomorphic balls in $\CC^k$, its Kobayashi metric identically vanishes and it is not biholomorphic to $\CC^k$. 

According to Theorem 1.12 in ~\cite{Fornaess}, $K_c=\{g^+<c\}$ for $c>0$ is a short $\CC^2$.
\begin{corollary}
The set $K_c$ with $c>0$ is a short $\CC^2$ with real analytic boundary foliated by injective Brody curves.
\end{corollary}
\end{remark}

\section{The case of finite compositions of generalized H\'enon mappings}
In this section, we explain how to apply the proof of Theorem \ref{thm:BrodyLeaf_intro} to finite compositions of generalized H\'enon mappings and give a sketch of the proof of Theorem \ref{thm:maingeneral}. The essential ingredients for the proof of Theorem \ref{thm:BrodyLeaf_intro} are Theorem \ref{thm:NonexistenceofHoloCurve} and computations in Section \ref{sec:analyticdiscs}.

We write $\mathfrak{f}=f_N\circ\cdots\circ f_1$ where $f_i(z, w)=(p_i(z)-a_i w, z)$, $p_i(z)$ is a monic polynomial of degree $d_i\geq 2$ and $a_i$ is a non-zero constant. By convention, we denote by $f_0=\mathrm{id}$. Let $d=\prod d_i$ and $a=\prod a_i$. The definitions of $\mathfrak{g}^+$, $K^\pm$, $U^\pm$ and $J^\pm$, and the extensions of ${f_i}^\pm$ to $\PP^k$ as a birational map are well-defined. (cf. \cite{DinhSibony}) Let $I_+=[0: 1: 0]$ and $I_-=[1: 0: 0]$. In this section, $K_c=\{\mathfrak{g}^+\leq c\}$.

\subsection{Theorem \ref{thm:NonexistenceofHoloCurve}} Let $K_{c, i}=f_i\circ \cdots \circ f_i(K_c)$. Since each ${f_i}^{\pm}$ is an automorphism of $\CC^2$ and  holomorphically extends to $\PP^k\setminus I_\pm\to \PP^k$ with $I_\mp$ as its fixed points such that ${f_i}^\pm(\{t=0\}\setminus I_\pm)=I_\mp$ (by abuse of notation, we used ${f_i}^\pm$ for their extension as well), respectively, $\overline{K_{c, i}}=K_{c, i}\cup I_+$.

Suppose to the contrary that $\phi:\Delta\to\PP^2$ be a non-trivial holomorphic curve such that $\phi(\Delta)\subset \overline{K_c}$ and $\phi(0)=I_+$. By the same argument as in Theorem \ref{thm:NonexistenceofHoloCurve}, we can form a sequence of non-trivial holomorphic curves $\{[z_{n, i}: 1: t_{n, i}]=f_i\circ\cdots\circ f_1\circ \mathfrak{f}^n\circ\phi:\Delta\to \overline{K_{d^nc, i}} \}$ where $n=0, 1, 2, \cdots$ and $i=0, \cdots, N-1$, and $z_{n, i}$ has vanishing order at least $d_{i+1}-1$ at $0$ for $i=0,\cdots, N-1$. By the same reasoning as in Theorem \ref{thm:NonexistenceofHoloCurve}, we can prove the same statement for $\mathfrak{f}$.

Hence, we obtain
\begin{theorem}\label{thm:NonexistenceofHoloCurve'}
Assume that $K_c$ is as above. There is no non-trivial holomorphic curve in $\PP^2$, which passes through $I_+$, and is supported in $\overline{K_c}\subseteq\PP^2$ for $c>0$.
\end{theorem}

\subsection{Computations in Section \ref{sec:analyticdiscs}} 
We first determine the constants used in the proof of Theorem \ref{thm:BrodyLeaf_intro}.

Note that due to Remark 1 after Proposition 2.2 in \cite{Favre2000}, we have Proposition \ref{prop:Favrecoordinate} for $\mathfrak{f}$ with $d=\prod d_i$ and $a=\prod a_i$. Hence, we can determine $r_\phi$ and $c_\phi$ in the same way. Concerning $c_{V_+}$ and $R>1$, we can find universal ones working for every $f_i$ with $i=1, \cdots, N$. We pick $c$ sufficiently large so that
\begin{itemize}
\item $c>\max_W g^++\sum_{i=1}^N\max_{f_N\circ\cdots\circ f_i(W)}g^+$;
\item $(|a|c_\phi/d+\max_{U_{I_-}}|r|+1)<c_\phi \exp(c)/8$ and $R<\exp(c)$.
\end{itemize}
The first condition is for our analytic discs not intersecting $W$ and the second condition is for the definition of our analytic discs.

We define as in the case of Theorem \ref{thm:BrodyLeaf_intro} a sequence of analytic discs
\begin{align*}
\phi_{s, n}={\mathfrak{f}}^{-n}\circ \psi^{-1}\circ\phi(s^{d^n}, c_\phi\theta/2+\theta_{s, n}),\quad\theta\in \Delta
\end{align*}
where $\psi, \phi$ are the functions in Proposition \ref{prop:Favrecoordinate}, and $\theta_{s, n}$ is defined in the same way as before.

Then, we can carry out the same estimates in Section \ref{sec:analyticdiscs} for ${f_i}^{-1}$ over $f_i\circ\cdots\circ f_1\circ \mathfrak{f}^m\circ\phi_{s, n}(\Delta)$ for $i=1, \cdots, N$ and for $0\leq m\leq n-1$. The common filtration $\{V^+, V^-, W\}$ makes it possible to do these estimates in the exactly same way as in the case of Theorem \ref{thm:BrodyLeaf_intro}.

For Lemma \ref{lem:V+} and also for other estimates in Section \ref{sec:analyticdiscs}, it is important to have $|z'|\leq \widetilde{c_g}|w'|$ for some fixed $\widetilde{c_g}>0$ where $(z, w)=f_i\circ\cdots\circ f_1\circ \mathfrak{f}^m\circ\phi_{s, n}(\theta)\in V^+$ and $(z', w')=[d(f_i\circ\cdots\circ f_1\circ \mathfrak{f}^m\circ\phi_{s, n})](\frac{d}{d\theta})$ at $(z, w)$. This is easy to see from Lemma \ref{lem:verticality} and 
\begin{align*}
D{f_i}^{-1}=\left(
\begin{array}{cc}
0 & 1\\ 1/a & {p_i}'(w)/a 
\end{array}\right)
\end{align*}
for $i=1, \cdots, N$ and for $w$ with $|w|>R$ as ${p_i}'(z)/a$ dominates other entries of the matrix.

Together with the above two subsections, the application of the same argument as in Section \ref{sec:proof} proves Theorem \ref{thm:maingeneral}. The density property is from \cite{Favre2000}.


\end{document}